\documentclass[preprint,12pt]{elsarticle}

\usepackage{amssymb}

\begin{document}

\begin{frontmatter}



\title{On the inverse problem of identifying  an unknown coefficient in a space-time fractional differential equation }


\author{Mine Aylin Bayrak$^{1,*}$}
\author{Ali Demir$^{1}$}
\address{$^{1}$Department of Mathematics, Kocaeli University,  Izmit, Kocaeli, Turkey}
\begin{abstract}
In this study, we focus on identifying solution and an unknown space-dependent coefficient in a space-time fractional differential equation by employing fractional Taylor series method. The substantial advantage of this method is that we don't take any over-measured data into account. Consequently, we determine the solution and unknown coefficient more precisely. The presented examples illustrate that outcomes of this method are in high agreement with the exact ones of the corresponding problem. Moreover, it can be implemented and applied effectively comparing with other methods.
\end{abstract}




\end{frontmatter}

\section{Introduction}
\label{}
Inverse problems of identifying unknown coefficients or source terms draw the growing attention of researchers in diverse branches of science, since it has diverse applications in real life \cite{1,2,3,4,5,6,7,8,9,10,11,12}. Identifying unknown parameters in fractional differential equations by numerical methods is one of the significant challenges in inverse problems \cite{13,14,15,16} since modelling scientific processes by fractional differential equations becomes one of the most attractive topics in science. One of the main reasons of this result is the nonlocal properties of fractional derivatives,  which leads to better modelling. \\
In this research, inverse problems of unknown coefficients in a space-time fractional differential equation with Neumann like boundary conditions are considered by utilizing fractional Taylor series method. The fractional derivative is in Caputo sense \cite{17}, which is one of the most common fractional derivatives. In the establishment of solution and space-dependent unknown coefficient, no need of any over-measured data is one of the significant aspect of this method, since using over-measured data leads to greater error in calculations. Moreover, the Neumann like boundary condition at the final point and initial condition allow us to establish the solution and unknown space-dependent coefficient without any difficulty. The other boundary and initial conditions  guarantee uniqueness of the solution and unknown coefficient. As a result, we conclude that having all these properties makes fractional Taylor series method one of the best method in inverse and direct problems of fractional differential equations. \\
The focus of this article is on the identifying of  the unknown space-dependent coefficient of the following governing  space-time fractional differential equation with initial and Neumann-like boundary conditions:
\begin{eqnarray}\label{1}
D_{t}^{\alpha}u(x,t)=D_{x}^{2\beta}u(x,t)+p(x)f(x,t),0<x<1,0<t<1, 0<\alpha,\beta\leq 1,
\end{eqnarray}
\begin{eqnarray}\label{2}
u(x,0)=\varphi(x),0\leqslant x\leqslant 1,
\end{eqnarray}
\begin{eqnarray}\label{3}
D_{x}^{\beta}u(0,t)=\mu_{1}(t),0< t\leqslant 1,
\end{eqnarray}
\begin{eqnarray}\label{4}
D_{x}^{\beta}u(1,t)=\mu_{2}(t),0< t\leqslant 1,
\end{eqnarray}
\section{Preliminaries}
Essential concepts and features of  fractional derivatives are presented in this section \cite{1,2,3,4}.\\
\textit{Definition 1.} The Riemann-Liouville fractional integral of order $\alpha~(\alpha\geqslant 0)$  is given as \\
\begin{eqnarray}\label{5}
J^{\alpha}f(x)=\frac{1}{\Gamma(\alpha)}\int_{0}^{x}(x-t)^{\alpha -1}f(t)dt,~\alpha>0,~x> 0,
\end{eqnarray}
\begin{eqnarray}\label{6}
J^{0}f(x)=f(x).
\end{eqnarray}
\textit{Definition 2.} The  Liouville-Caputo fractional derivative of order $\alpha$  is given as  \\
\begin{eqnarray}\label{7}
D^{\alpha}f(x)=J^{n-\alpha}D^{n}f(x)=\int_{0}^{x}(x-t)^{n-\alpha -1}\frac{d^{n}}{dt^{n}}f(t)dt,~n-1<\alpha<n,~x> 0,
\end{eqnarray}
where $D^{n}$ denotes the ordinary derivative of order $n$.\\
\textit{Definition 3.} The $\alpha^{th}$ order derivative of $u(x,t)$ in Liouville-Caputo sense is given as  \\
\begin{eqnarray}\label{8}
D_{t}^{\alpha}u(x,t)=
\left\{
 \begin{array}{lr}
\frac{1}{\Gamma(n-\alpha)}\int_{0}^{t}(t-\xi)^{n-\alpha -1}\frac{\partial^{n}u(x,\xi)}{\partial t^{n}}d\xi,~n-1<\alpha<n,\\
\frac{\partial^{n}u(x,t)}{\partial t^{n}}~~~~~~~~~~~~~~~~~~~~~~~~~~~~~~~~~~,\alpha=n\in N.
\end{array}
   \right.
\end{eqnarray}
\textit{Definition 4.} An $(\alpha,\beta)$-fractional Taylor series is defined as follows \cite{18}:
\begin{eqnarray}\label{9}
\sum\limits_{i+j=0}^{\infty}g_{i,j}t^{i\alpha}x^{j\beta}=\underbrace{g_{0,0}}_{i+j=0}+\underbrace{g_{1,0}t^{\alpha}+g_{0,1}x^{\beta}}_{i+j=1}+
...+\underbrace{\sum\limits_{k=0}^{n}g_{n-k,k}t^{(n-k)\alpha}x^{k\beta}}_{i+j=n}+...
\end{eqnarray}
where $g_{i,j}, i,j\epsilon N$ are the coefficients of the series.\\
\textit{Lemma 5.} Let $u(x,t)$ has a fractional Taylor series  representation as (\ref{9}) for $(x,t)\epsilon [0,R_{x})\times [0,R_{t})$. If $D_{t}^{r\alpha}D_{x}^{s\beta}u(x,t)\epsilon ((0,R_{x})\times (0,R_{t}))$ for $r,s\epsilon N$, then
\begin{eqnarray}\label{10}
D_{t}^{r\alpha}u(x,t)=\sum\limits_{i+j=0}^{\infty}g_{i+r,j}\frac{\Gamma((i+r)\alpha+1)}{\Gamma(i\alpha+1)}t^{i\alpha}x^{j\beta}
\end{eqnarray}
\begin{eqnarray}\label{11}
D_{x}^{s\beta}u(x,t)=\sum\limits_{i+j=0}^{\infty}g_{i,j+s}\frac{\Gamma((j+s)\beta+1)}{\Gamma(j\beta+1)}t^{i\alpha}x^{j\beta}
\end{eqnarray}
\section{Fractional Taylor series method}
In order to determine the  diffusion coefficient $p(x)$ of time in the space-time fractional diffusion problem (\ref{1})-(\ref{4}), in the series form we plug the  fractional Taylor series of $u=u(x,t)$ and $p=p(x)$ into (\ref{1})-(\ref{4}) which leads to:
\begin{eqnarray}\label{12}
\sum\limits_{i+j=0}^{\infty}g_{i+1,j}\frac{\Gamma((i+1)\alpha+1)}{\Gamma(i\alpha+1)}t^{i\alpha}x^{j\beta}
=\sum\limits_{i+j=0}^{\infty}g_{i,j+2}\frac{\Gamma((j+2)\beta+1)}{\Gamma(j\beta+1)}t^{i\alpha}x^{j\beta}\nonumber\\
+\sum\limits_{k=0}^{\infty}p_{k}\Big\{\sum\limits_{i+j=0}^{\infty}\frac{t^{i\alpha}x^{(j+k)\beta}}{\Gamma(i\alpha+1)\Gamma((j+k)\beta+1)}\Big\}
\end{eqnarray}
Making two series on both sides of above equation equal to each other, the unknown coefficients in the fractional Taylor series  of $p(x)$ are acquired.\\
\section{Illustrative Examples}
In this section, we present some examples to illustrate implementation of fractional Taylor series method for the inverse problem of revealing an unknown space-dependent coefficient in space-time fractional differential equations. The following examples are about determination of unknown coefficient depending on $x$.\\
\textit{Example 1.}  Consider the inverse space-dependent coefficient problem involving space-time fractional differential equations:
\begin{eqnarray}\label{13}
D_{t}^{\alpha}u(x,t)=D_{x}^{2\beta}u(x,t)+p(x)u(x,t),0<x<1,0<t<1, 0<\alpha,\beta\leq 1,
\end{eqnarray}
\begin{eqnarray}\label{14}
u(x,0)=E_{\beta}(x^{2\beta}),0\leqslant x\leqslant 1,
\end{eqnarray}
\begin{eqnarray}\label{15}
D_{x}^{\beta}u(0,t)=0,0< t\leqslant 1,
\end{eqnarray}
\begin{eqnarray}\label{16}
D_{x}^{\beta}u(1,t)=E_{\alpha}(2t^{\alpha})D_{x}^{\beta}E_{\beta}(x^{2\beta})|_{x=1},0< t\leqslant 1,
\end{eqnarray}
where $E_{\beta}(x^{2\beta})=\sum\limits_{j=1}^{\infty}\frac{x^{2j\beta}}{\Gamma(j\beta+1)}$ is the fractional generalization of the function $exp(x^{2})$. We determine the unknown function $p(x)$ in fractional Taylor series form as follows:
\begin{eqnarray}\label{17}
p(x)=\sum\limits_{k=0}^{\infty}p_{k}\frac{x^{k\beta}}{\Gamma(1+k\beta)}, 0<\beta\leq 1.
\end{eqnarray}
which leads to Eq. (\ref{12}),with the initial coefficients
\begin{eqnarray}\label{18}
g_{0,j}=\frac{1}{\Gamma(j\beta+1)},
\end{eqnarray}
\begin{eqnarray}\label{19}
g_{0,0}=1.
\end{eqnarray}
The coefficients $g_{i,j}$ are acquired by equating two series in Eq. (\ref{12}), which
allow us to form the solution of Eq.(\ref{13}) as follows:
\begin{eqnarray}\label{20}
u(x,t)=1+\frac{\Gamma(2\beta+1)}{\Gamma(\beta+1)}\frac{t^{\alpha}}{\Gamma(\alpha+1)}+\frac{\Gamma(2\beta+1)}{\Gamma(\beta+1)}\frac{t^{2\alpha}}{\Gamma(2\alpha+1)}\nonumber\\
+4\frac{\Gamma(2\beta+1)}{\Gamma(\beta+1)}\frac{t^{3\alpha}}{\Gamma(3\alpha+1)}+2\frac{\Gamma(2\beta+1)}{\Gamma(\beta+1)}\frac{t^{\alpha}}{\Gamma(\alpha+1)}\frac{x^{2\beta}}{\Gamma(2\beta+1)}\nonumber\\
+8\frac{\Gamma(2\beta+1)}{\Gamma(\beta+1)}\frac{t^{4\alpha}}{\Gamma(4\alpha+1)}+4\frac{\Gamma(2\beta+1)}{\Gamma(\beta+1)}\frac{t^{2\alpha}}{\Gamma(2\alpha+1)}\frac{x^{2\beta}}{\Gamma(2\beta+1)}\nonumber\\
+\frac{\Gamma(4\beta+1)}{\Gamma(2\beta+1)}\frac{x^{4\beta}}{\Gamma(4\beta+1)}+16\frac{\Gamma(2\beta+1)}{\Gamma(\beta+1)}\frac{t^{5\alpha}}{\Gamma(5\alpha+1)}\nonumber\\
+8\frac{\Gamma(2\beta+1)}{\Gamma(\beta+1)}\frac{t^{3\alpha}}{\Gamma(3\alpha+1)}\frac{x^{2\beta}}{\Gamma(2\beta+1)}+32\frac{\Gamma(2\beta+1)}{\Gamma(\beta+1)}\frac{t^{6\alpha}}{\Gamma(6\alpha+1)}\nonumber\\
+16\frac{\Gamma(2\beta+1)}{\Gamma(\beta+1)}\frac{t^{4\alpha}}{\Gamma(4\alpha+1)}\frac{x^{2\beta}}{\Gamma(2\beta+1)}+4\frac{\Gamma(2\beta+1)}{\Gamma(\beta+1)}\frac{t^{2\alpha}}{\Gamma(2\alpha+1)}\frac{x^{4\beta}}{\Gamma(4\beta+1)}\nonumber\\
+\frac{\Gamma(2\beta+1)}{\Gamma(\beta+1)}\frac{x^{6\beta}}{\Gamma(6\beta+1)}
\end{eqnarray}
In order to establish the unknown coefficient $p(x)$, the Neumann-like  boundary condition
at $x = 1$ taken into account in (\ref{16}) which produces the coefficients $p_{k}$ as follows:\\
$p_{0}=0$,\\
$p_{1}=0$,\\
$p_{2}=2\frac{\Gamma(2\beta+1)}{\Gamma(\beta+1)}-\frac{\Gamma(4\beta+1)}{\Gamma(2\beta+1)}$,\\
\vdots\\
As a result, the unknown  coefficient $p(x)$ is determined in the series form as follows:
\begin{eqnarray}\label{21}
p(x)=\Big(2\frac{\Gamma(2\beta+1)}{\Gamma(\beta+1)}-\frac{\Gamma(4\beta+1)}{\Gamma(2\beta+1)}\Big)\frac{x^{2\beta}}{\Gamma(2\beta+1)}+...
\end{eqnarray}
\newpage
\setlength{\tabcolsep}{0.009em} 
{\renewcommand{\arraystretch}{1.3}
\begin{table}
  \caption{Comparison of absolute errors  at $x=0.5$   with $E(\alpha,\beta)$ of Example 1. }
  \begin{center}
  \begin{tabular}{c|c|c|c|c|c|c|c|c|c|c}
	 \hline
 $t$&$Exact$ &$E(1,1)$&$E(1,0.9)$&$E(1,0.7)$&$E(0.9,1)$&$E(0.9,0.9)$&$E(0.9,0.7)$&$E(0.7,1)$&$E(0.7,0.9)$&$E(0.7,0.7)$\\
 \hline
    0.05 & 	 1.41888 &	4.40e-02 &	4.14e-02 &	4.16e-02 &	4.39e-02 &	3.51e-02 &	2.57e-02 &	4.32e-02 &	1.13e-02 &	3.48e-02\\
    0.10  &	 1.56810 &	4.36e-02 &	2.57e-02 &	1.84e-03 &	4.33e-02 &	1.46e-02 &	2.66e-02 &	4.13e-02 &	2.74e-02 &	1.33e-01\\
    0.15  &	 1.73302 &	4.32e-02 &	8.34e-03 &	4.24e-02 &	4.26e-02 &	7.99e-03 &	8.38e-02 &	3.81e-02 &	7.06e-02 &	2.41e-01\\
    0.20  &	 1.91528 &	4.26e-02 &	1.11e-02 &	9.17e-02 &	4.16e-02 &	3.31e-02 &	1.48e-01 &	3.26e-02 &	1.21e-01 &	3.66e-01\\
    0.25  &	 2.11672 &	4.18e-02 &	3.27e-02 &	1.47e-01 &	4.00e-02 &	6.14e-02 &	2.19e-01 &	2.41e-02 &	1.81e-01 &	5.11e-01\\
    0.30  &	 2.33933 &	4.06e-02 &	5.70e-02 &	2.08e-01 &	3.78e-02 &	9.35e-02 &	3.00e-01 &	1.14e-02 &	2.52e-01 &	6.82e-01\\
    0.35  &	 2.58534 &	3.91e-02 &	8.43e-02 &	2.77e-01 &	3.47e-02 &	1.30e-01 &	3.91e-01 &	6.69e-03 &	3.37e-01 &	8.83e-01\\
    0.40  &	 2.85721 &	3.70e-02 &	1.15e-01 &	3.54e-01 &	3.04e-02 &	1.72e-01 &	4.94e-01 &	3.14e-02 &	4.39e-01 &	1.12e+00\\
    0.45  &	 3.15763 &	3.42e-02 &	1.50e-01 &	4.41e-01 &	2.46e-02 &	2.19e-01 &	6.12e-01 &	6.44e-02 &	5.60e-01 &	1.40e+00\\
    0.50  &	 3.48959 &	3.06e-02 &	1.89e-01 &	5.38e-01 &	1.68e-02 &	2.74e-01 &	7.46e-01 &	1.07e-01 &	7.04e-01 &	1.72e+00\\
\hline
\end{tabular}
\end{center}
\end{table}
\textit{Example 2.}  Consider the inverse coefficient problem involving space-time fractional differential equations:
\begin{eqnarray}\label{22}
D_{t}^{\alpha}u(x,t)=D_{x}^{2\beta}u(x,t)+p(x)u(x,t),0<x<1,0<t<1, 0<\alpha,\beta\leq 1,
\end{eqnarray}
\begin{eqnarray}\label{23}
u(x,0)=E_{\beta}(x^{3\beta}),0\leqslant x\leqslant 1,
\end{eqnarray}
\begin{eqnarray}\label{24}
D_{x}^{\beta}u(0,t)=0,0< t\leqslant 1,
\end{eqnarray}
\begin{eqnarray}\label{25}
D_{x}^{\beta}u(1,t)=E_{\alpha}(t^{\alpha})D_{x}^{\beta}E_{\beta}(x^{3\beta})|_{x=1},0< t\leqslant 1,
\end{eqnarray}
where $E_{\beta}(x^{3\beta})=\sum\limits_{j=1}^{\infty}\frac{x^{3j\beta}}{\Gamma(j\beta+1)}$ is the fractional generalization of the function $exp(x^{3})$. We determine the unknown function $p(x)$ in fractional Taylor series form as follows:
\begin{eqnarray}\label{26}
p(x)=\sum\limits_{k=0}^{\infty}p_{k}\frac{x^{k\beta}}{\Gamma(1+k\beta)}, 0<\beta\leq 1.
\end{eqnarray}
which leads to Eq. (\ref{12}),
with the initial coefficients
\begin{eqnarray}\label{27}
g_{0,j}=\frac{1}{\Gamma(j\beta+1)},
\end{eqnarray}
\begin{eqnarray}\label{28}
g_{0,0}=1.
\end{eqnarray}
The coefficients $g_{i,j}$ are acquired by equating two series in Eq. (\ref{12}), which
allow us to form the solution of Eq.(\ref{22}) as follows:
\begin{eqnarray}\label{29}
u(x,t)=1+(1+\frac{\Gamma(2\beta+1)}{\Gamma(\beta+1)})\frac{t^{\alpha}}{\Gamma(\alpha+1)}+(1+\frac{\Gamma(2\beta+1)}{\Gamma(\beta+1)})\frac{t^{2\alpha}}{\Gamma(2\alpha+1)}\nonumber\\
+(1+\frac{\Gamma(2\beta+1)}{\Gamma(\beta+1)})\frac{t^{3\alpha}}{\Gamma(3\alpha+1)}+\frac{\Gamma(3\beta+1)}{\Gamma(\beta+1)}\frac{x^{3\beta}}{\Gamma(3\beta+1)}\nonumber\\
+(1+\frac{\Gamma(2\beta+1)}{\Gamma(\beta+1)})\frac{t^{4\alpha}}{\Gamma(4\alpha+1)}+\frac{\Gamma(3\beta+1)}{\Gamma(\beta+1)}\frac{t^{\alpha}}{\Gamma(\alpha+1)}\frac{x^{3\beta}}{\Gamma(3\beta+1)}\nonumber\\
+(1+\frac{\Gamma(2\beta+1)}{\Gamma(\beta+1)})\frac{t^{5\alpha}}{\Gamma(5\alpha+1)}+\frac{\Gamma(3\beta+1)}{\Gamma(\beta+1)}\frac{t^{2\alpha}}{\Gamma(2\alpha+1)}\frac{x^{3\beta}}{\Gamma(3\beta+1)}\nonumber\\
+(1+\frac{\Gamma(2\beta+1)}{\Gamma(\beta+1)})\frac{t^{6\alpha}}{\Gamma(6\alpha+1)}+\frac{\Gamma(3\beta+1)}{\Gamma(\beta+1)}\frac{t^{3\alpha}}{\Gamma(3\alpha+1)}\frac{x^{3\beta}}{\Gamma(3\beta+1)}\nonumber\\
+\frac{\Gamma(6\beta+1)}{\Gamma(2\beta+1)}\frac{x^{6\beta}}{\Gamma(6\beta+1)}
\end{eqnarray}
In order to establish the unknown coefficient $p(x)$, the Neumann boundary condition
at $x = 1$ taken into account in (\ref{25}) which produces the coefficients $p_{k}$ as follows:\\
$p_{0}=1$,\\
$p_{1}=-\frac{\Gamma(3\beta+1)}{\Gamma(\beta+1)}$,\\
$p_{2}=0$,\\
$p_{3}=-\frac{\Gamma(6\beta+1)}{\Gamma(2\beta+1)}+\frac{\Gamma(3\beta+1)\Gamma(4\beta+1)}{(\Gamma(\beta+1))^{2}}$,\\
\vdots\\
As a result, the unknown  coefficient $p(x)$ is determined in the series form as follows:
\begin{eqnarray}\label{30}
p(x)=1-\frac{\Gamma(3\beta+1)}{\Gamma(\beta+1)}\frac{x^{\beta}}{\Gamma(\beta+1)}\nonumber\\
-\Big(\frac{\Gamma(6\beta+1)}{\Gamma(2\beta+1)}-\frac{\Gamma(3\beta+1)\Gamma(4\beta+1)}{(\Gamma(\beta+1))^{2}}\Big)\frac{x^{2\beta}}{\Gamma(2\beta+1)}+...
\end{eqnarray}

\newpage
\setlength{\tabcolsep}{0.009em} 
{\renewcommand{\arraystretch}{1.3}
\begin{table}
  \caption{Comparison of absolute errors  at $x=1$   with $E(\alpha,\beta)$ of Example 2. }
  \begin{center}
  \begin{tabular}{c|c|c|c|c|c|c|c|c|c|c}
	 \hline
 $t$&$Exact$ &$E(1,1)$&$E(1,0.9)$&$E(1,0.7)$&$E(0.9,1)$&$E(0.9,0.9)$&$E(0.9,0.7)$&$E(0.7,1)$&$E(0.7,0.9)$&$E(0.7,0.7)$\\
 \hline
   0.005  	&  2.51253  &	 7.52e-03 &	 5.75e-03 &	 2.82e-03 &	 1.33e-02 &	 1.02e-02 &	 4.99e-03 &	 4.12e-02 &	 3.15e-02 &	 1.54e-02\\
   0.010  	&  2.52513  &	 1.51e-02 &	 1.15e-02 &	 5.65e-03 &	 2.49e-02 &	 1.91e-02 &	 9.35e-03 &	 6.77e-02 &	 5.17e-02 &	 2.54e-02\\
   0.015  	&  2.53778  &	 2.27e-02 &	 1.73e-02 &	 8.49e-03 &	 3.61e-02 &	 2.76e-02 &	 1.35e-02 &	 9.08e-02 &	 6.94e-02 &	 3.40e-02\\
   0.020  	&  2.55050  &	 3.03e-02 &	 2.32e-02 &	 1.14e-02 &	 4.69e-02 &	 3.59e-02 &	 1.76e-02 &	 1.12e-01 &	 8.56e-02 &	 4.20e-02\\
   0.025  	&  2.56329  &	 3.80e-02 &	 2.90e-02 &	 1.42e-02 &	 5.76e-02 &	 4.40e-02 &	 2.16e-02 &	 1.32e-01 &	 1.01e-01 &	 4.95e-02\\
   0.030  	&  2.57614  &	 4.57e-02 &	 3.49e-02 &	 1.71e-02 &	 6.81e-02 &	 5.21e-02 &	 2.55e-02 &	 1.51e-01 &	 1.16e-01 &	 5.66e-02\\
   0.035  	&  2.58905  &	 5.34e-02 &	 4.08e-02 &	 2.00e-02 &	 7.85e-02 &	 6.00e-02 &	 2.94e-02 &	 1.70e-01 &	 1.30e-01 &	 6.35e-02\\
   0.040  	&  2.60203  &	 6.12e-02 &	 4.68e-02 &	 2.29e-02 &	 8.89e-02 &	 6.79e-02 &	 3.33e-02 &	 1.88e-01 &	 1.43e-01 &	 7.03e-02\\
   0.045  	&  2.61507  &	 6.90e-02 &	 5.28e-02 &	 2.59e-02 &	 9.92e-02 &	 7.58e-02 &	 3.71e-02 &	 2.05e-01 &	 1.57e-01 &	 7.68e-02\\
   0.050  	&  2.62818  &	 7.69e-02 &	 5.88e-02 &	 2.88e-02 &	 1.09e-01 &	 8.36e-02 &	 4.10e-02 &	 2.22e-01 &	 1.70e-01 &	 8.32e-02\\
\hline
\end{tabular}
\end{center}
\end{table}

\newpage
\section{Conclusion}
In this research, the identification of an unknown space-time coefficient  in a space-time fractional differential equation with initial and Neumann-like boundary conditions by fractional Taylor series method is investigated. First, the solution of the equation is established by employing fractional Taylor series method with the help of initial condition and Neumann-like boundary condition at the final point. Later, the unknown space-dependent coefficient is identified in the series form without using any over-measured data which is substantial advantage of this method. Future work will be on the construction of unknown coefficient or source term in space-dependent fractional differential equations with various boundary conditions.

\section{Acknowledgements}
The authors are grateful to the referees for their comments and suggestions which have improved the paper. This research is not funded by any organization.




\newpage
\textbf{References}
\begin{thebibliography}{99}

\bibitem{1} Kilbas, A.A. , Srivastava, H.M. and Trujillo, J.J. Theory and applications of fractional  differential equations, Amsterdam: Elsevier, 2006.
\bibitem{2} Podlubny, I. Fractional differential equation, San Diego, CA: Academic Press, 1999.
\bibitem{3} Sabatier, J., Agarwal, O.P. and Machado, J.A.T.(eds). Advances in fractional calculus: theoretical developments and applications in physics and  engineering, Dordrecht: Springer, 2007.
\bibitem{4} Samko, S.G., Kilbas, A.A. and Marichev, O.I. Fractional integrals and derivatives theory and applications, Amsterdam: Gordon and Breach, 1993.
 \bibitem{5} Odibat, Z. Approximations of fractional integrals and Caputo fractional derivatives. Appl. Math. Comput.178, 527-533, 2006.
\bibitem{6} Momani, S., Odibat, Z. Analytical approach to linear fractional partial differential equations arising in fluid mechanics. Phys. Lett. A 355, 271-279, 2006.
\bibitem{7} Seki, K., Wojcik, M., Tachiya, M. Fractional reaction-diffusion equation. J. Chem. Phys. 119, 2165, 2003.
\bibitem{8} Saxena, R.K., Mathai, A.M., Haubold, H.J. Fractional reaction–diffusion equations. Astrophys. Space Sci. 305, 289-296, 2006.
\bibitem{9} Abbasbandy, S. The application of homotopy analysis method to nonlinear equations arising in heat transfer, Physics Lett. A  360, 109-113, 2006.
\bibitem{10} Atangana, A. On the new fractional derivative and application to nonlinear Fisher-reaction- diffusion equation, Appl. Math. Comput. 273, 948-956, 2016.
 \bibitem{11} El-Ajou, A., Abu Arqub, O., Momani, S., Baleanu, D., Alsaedi, A. A novel expansion iterative method for solving linear partial differential equations of fractional order. Appl. Math. Comput. 257, 119-133, 2015.
\bibitem{12} Bayrak, M.A., Demir, A. A new approach for space-time fractional partial differential equations by residual power series method. Appl. Math. Comput. 336, 215-230, 2018.
\bibitem{13}Xiangtuan X. , Hongbo G. ,Xiaohong Liu., An inverse problem for a fractional diffusion equation,Journal of Computational and Applied Mathematics 236, 4474–4484, 2012.
\bibitem{14}Mansur I. I. , Muhammed C., Inverse source problem for a time-fractional diffusion equation with nonlocal boundary conditions, Applied Mathematical Modelling 40, 4 891-899, 2016.
\bibitem{15} Songshu L.  and Lixin F. ,   An Inverse Problem for a Two-Dimensional Time-Fractional Sideways Heat Equation, Mathematical Problems in Engineering, Volume 2020, Article ID 5865971, 13 pages.
\bibitem{16}Zhiyuan L., Xing C., and Gongsheng L., An inverse problem in time-fractional diffusion equations with nonlinear boundary condition, Journal of Mathematical Physics, 60, 091502 ,2019.
\bibitem{17}Jaradat I, Alquran M. , Abdel-Muhsen R., An analytical framework of 2D diffusion, wave-like, telegraph and Burgers models with twofold Caputo derivatives ordering, Nonlinear Dynamics, 93(4) 1911-1922, 2018.
\bibitem{18} Abdel-Muhsen, R., An application of Taylor series method in higher dimensional fractal spaces, The 6th International Arap conference on Mathematics and Computation, 2019.




\end{thebibliography}


\end{document}